\documentclass[letterpaper, 10 pt, conference]{ieeeconf}  

\IEEEoverridecommandlockouts                              

\usepackage{graphics} 
\usepackage{epsfig} 

\usepackage{amsmath, amsthm} 
\usepackage{amssymb}  
\usepackage[mathscr]{euscript}   
\usepackage[usenames,dvipsnames]{color} 
\usepackage{algpseudocode}
\usepackage{algorithm}
\usepackage{hyperref}
\usepackage{inconsolata}

\newtheoremstyle{main}
    {1em}                                                
    {1em}                                                
    {\itshape}                                        
    {0pt}                                                
    {\scshape}                                           
    {\\*}                                                
    {2pt}                                                
    {\thmname{#1}\thmnumber{ #2}: \thmnote{\itshape #3}}

\theoremstyle{main}
\newtheorem{theorem}{Theorem}
\newtheorem{proposition}{Proposition}
\newtheorem{corollary}{Corollary}
\theoremstyle{definition} 
\newtheorem{definition}{Definition}

\newcommand\eqdef{\mathrel{\overset{\makebox[0pt]{\mbox{\normalfont\tiny\sffamily def}}}{=}}}

\title{\LARGE \bf
Combinatorial Disjunctive Constraints for Obstacle Avoidance in Path Planning
}

\author{Raul Garcia$^{1}$, Illya V. Hicks$^{2}$ and Joey Huchette$^{3}$
\thanks{$^{1}$Raul Garcia is with the Department of Computational Applied Mathematics \& Operations Research, Rice 
        University, Houston, TX 77005, USA
        {\tt\small rjgarcia@rice.edu}}%
\thanks{$^{2}$Illya V. Hicks is with Faculty of Computational Applied Mathematics \& Operations Research, Rice      
        University, Houston, TX 77005, USA
        {\tt\small ivhicks@rice.edu}}%
\thanks{$^{3}$Joey Huchette is with Google Research,
        Cambridge, MA 02142, USA
        {\tt\small jhuchette@google.com}}%
}

\begin{document}

\maketitle
\thispagestyle{empty}
\pagestyle{empty}

\begin{abstract}

We present a new approach for modeling avoidance constraints in 2D environments, in which waypoints are assigned to obstacle-free polyhedral regions. Constraints of this form are often formulated as mixed-integer programming (MIP) problems employing big-M techniques---however, these are generally not the strongest formulations possible with respect to the MIP's convex relaxation (so called \emph{ideal} formulations), potentially resulting in larger computational burden. We instead model obstacle avoidance as combinatorial disjunctive constraints and leverage the \emph{independent branching} scheme to construct small, ideal formulations. As our approach requires a biclique cover for an associated graph, we exploit the structure of this class of graphs to develop a fast subroutine for obtaining biclique covers in polynomial time. We also contribute an open-source Julia library named \texttt{ClutteredEnvPathOpt} to facilitate computational experiments of MIP formulations for obstacle avoidance. Experiments have shown our formulation is more compact and remains competitive on a number of instances compared with standard big-M techniques, for which solvers possess highly optimized procedures.

\end{abstract}

\section{INTRODUCTION}


An essential task in many robotics applications is for an autonomous agent to find a path through a cluttered environment to safely arrive at a destination. Paths are often represented as a sequence of waypoints satisfying dynamic and obstacle avoidance constraints, and may or may not integrate time. Approaches to path planning can be broadly categorized into discrete searches on graphs or continuous optimization.

Discrete approaches typically utilize graphs to capture feasible paths among a set of waypoints within the environment \cite{Lozano-Perez, Masehian} or may construct a tree of possible paths given a successor set of available actions \cite{Baudouin2011, Kuffner2001, Chestnutt2003, Michel2005}. Shortest paths between the start and goal positions can then computed via Dijkstra's algorithm or $A^*$. Obstacle avoidance is accomplished by appropriately prohibiting certain edges/nodes in the former or pruning the tree accordingly in the latter. While these approaches have worked well in certain instances, they tend to suffer from computational burden in complex environments and zig-zagging behavior \cite{Masehian, Baudouin2011, SurveyMotionPlanning}, and they do not consider the optimality of the exact waypoint locations. Sampling-based algorithms such as Probabilistic Roadmaps \cite{PRMs} and Rapidly-exploring Random Trees \cite{RRT1998} help relieve computational burden; however, their probabilistic completeness property implies sampling may need to grow large before a feasible path may be found.

On the other hand, continuous optimization avoids the challenges of exploring graphs and instead directly operates on the waypoint placements via continuous decision variables. This allows the construction of a more expressive objective function, and dynamic constraints can be employed to prevent infeasible or undesirable behavior. However, incorporating agent yaw and obstacle avoidance generally results in nonlinearity and/or nonconvexity, presenting challenges in computational efficiency, global optimality, and outright guaranteeing obstacle avoidance \cite{Fallon, Herdt}.

The past two decades have seen a growing interest in mixed-integer programming (MIP) in the path planning community due to its ability to model nonlinearity and nonconvexity \cite{SurveyMotionPlanning, RichardsHow2002, Schouwenaars, MIPformulations}, in conjunction with the availability of specialized solvers. This is particularly useful for modeling collision avoidance, in which we require an agent to lie in some region of obstacle-free safe terrain, typically nonconvex due to cavities brought forth by the obstacles. 
One approach to modeling obstacle avoidance is the direct assignment of waypoints (or trajectory pieces) to obstacle-free polyhedral regions in the configuration space, through the use of binary variables \cite{FootstepDeits, DeitsUAV}. That is, letting $x$ represent a waypoint, we require
\begin{equation}
    \label{eq:general disjunction}
    x \in \bigcup_{i=1}^d P^i ,
\end{equation}
where each $P^i \subseteq \mathbb{R}^n$ is a polyhedron. These are \emph{disjunctive constraints} and can be modeled for each waypoint using MIP, potentially along other dynamic and/or logical constraints. Modeling this nonconvex domain has computational implications, however. The typical approach to formulating disjunctive constraints is through a big-M technique \cite{MIPformulations, FootstepDeits}. While these contraints are simple to implement and reason about, they generally have weak convex relaxations, which can often result in poor computational performance.

In this work, we leverage the \emph{independent branching} (IB) scheme introduced by Vielma and Nemhauser \cite{IBoriginal} and further expanded by Huchette and Vielma \cite{CDC}, in order to construct small, ideal formulations of combinatorial disjunctive constraints (CDC) for obstacle avoidance in 2D environments. We discuss conditions for when this approach may be employed, and, as this method requires a \emph{biclique cover} of an associated graph, we develop a polynomial time algorithm for finding biclique covers on this class of graphs. Furthermore, we contribute an open-source Julia library named \texttt{ClutteredEnvPathOpt}\footnote{\url{https://github.com/raulgarcia66/ClutteredEnvPathOpt.jl}} to facilitate computational experiments of MIP formulations for obstacle avoidance and
demonstrate our framework through an application of footstep planning for a humanoid robot, largely adapted from the mixed-integer quadratically constrained quadratic program
(MIQCQP) 
presented by Deits and Tedrake \cite{FootstepDeits}. While our experiments have shown the big-M method to generally outperform our IB scheme approach, we hypothesize it is due to the specialized heuristics of MIP solvers. Moreover, we highlight that in a handful of instances our method was able to compute optimal paths 40-50\% faster compared to big-M and remains competitive for many more.  
We hope our approach may become a favorable alternative in particular classes of instances.

This paper is organized as follows. In Section \ref{sec:prelim}, we provide background on MIP formulation theory and introduce the big-M and IB scheme approaches for disjunctive constraints. In Section \ref{sec:contributions}, we discuss conditions for when the IB scheme may be employed in 2D environments and present our algorithm for obtaining biclique covers on the class of graphs arising in this context. We present our computational results and an example of our framework in Section \ref{sec:results}, before ending with a discussion on the observed performance.

\section{PRELIMINARIES}
\label{sec:prelim}

\subsection{Polyhedral Theory}

A (bounded) $\mathscr{V}$-polyhedron 
is a set $P \subset \mathbb{R}^n$ that can be expressed as 
$$
P = \text{Conv} (V) \eqdef \left\{ \sum\limits_{v \in V} \lambda_v v : \lambda \in \Delta^{V} \right\}
$$
for some finite set of vectors $V \subset \mathbb{R}^n$, where $\Delta^{V} \eqdef \left\{ \lambda \in \mathbb{R}_{+}^{|V|} : \sum_{v \in V} \lambda_v = 1 \right\}$ is the standard simplex. In other words, it is the convex hull of the vectors $V$. We will commonly refer to bounded polyhedrons as \emph{polytopes}.

A $\mathscr{H}$-polyhedron 
is a set $P \subset \mathbb{R}^n$ that can be expressed as
$$
P = \{ x \in \mathbb{R}^n :  A x \leq b \}
$$
where $A$ is an $m \times n$ matrix and $b$ is a vector in $\mathbb{R}^m$. Without loss of generality, bounded polyhedra have both $\mathscr{V}$- and $\mathscr{H}$- representations. 

A (binary) MIP formulation $F$ for a constraint $x \in Q \subseteq \mathbb{R}^{n_1}$ is the composition of linear inequalities
\begin{equation}
    \begin{split}
    R \eqdef \{ (x,y,z) \in \mathbb{R}^{n_1} \times \mathbb{R}^{n_2} \times [0,1]^{n_3} : \\
    Ax + By + Cz \leq d \}
    \end{split}
\end{equation}
with (binary) integrality conditions
\begin{equation}
    F \eqdef R \cap (\mathbb{R}^{n_1} \times \mathbb{R}^{n_2} \times \{0,1\}^{n_3})
\end{equation}
such that $\text{Proj}_x (F) = Q$. The variables $y$ and $z$ denote auxiliary continuous and binary variables, respectively, which are employed to handle the nonlinearity and/or nonconvexity of $Q$. $R$ is referred to as the \emph{LP relaxation}, or just \emph{relaxation}, of $F$, since it relaxes the integrality conditions. 

For computational purposes, we are interested in understanding both the strength of a given formulation, as well as in quantifying its size or complexity.
A formulation is said to be \emph{ideal} if each extreme point of the relaxation naturally satisfies the integrality conditions.
This is the strongest possible MIP formulation with respect to the convex relaxation that one can expect; in particular, if an MIP formulation is ideal, it can be solved as a linear program (LP), foregoing the need for branching on fractional variables. 
On the other hand, one measure of the complexity of a formulation is the number of auxiliary continuous variables $y$ and auxiliary binary variables $z$ used, as well as the number of inequalities in the description of $R$. 




\subsection{Big-M Techniques}
\label{subsec:bigM}

A standard approach to formulating polyhedral disjunctive constraints is the big-M method \cite{MIPformulations}. 
Letting $P^i = \{ x \in \mathbb{R}^n: A^i x \leq b^i \}$, the disjunctive constraint \eqref{eq:general disjunction} can be modeled by
\begin{subequations}
    \label{eq:bigM formulation}
    \begin{align}
    & \quad A^i x \leq b^i z_i + M^i (1-z_i) & \forall i \in [d] \label{subeq:bigM} \\
    & \quad \sum_{i=1}^d z_i = 1 \\
    & \quad z \in \{0,1\}^d ,
    \end{align}
\end{subequations}
where $[d] \eqdef \{1,\ldots,d \}$ and $M^i, i \in [d],$ are vectors with sufficiently large entries. In particular, for each $i \in [d]$, we must take
\begin{equation}
    \label{eq:lower bound on M}
    M^i_k \geq \max\limits_{x \in \Omega} (A^i x)_k ,
\end{equation}
for each row index $k$
of the system $A^i x \leq b^i$ and $\Omega \eqdef \bigcup_{i=1}^d P^i$. $M$ values can be computed efficiently via linear programming by maximizing over some polytope containing $\Omega$.

Such formulations are not ideal in general; in fact, the values of $M$ play a crucial role in the strength of a formulation \cite{MIPformulations}. Larger values increase the size of the convex relaxation, which generally results in a larger branch-and-cut search tree. Therefore, it is preferable to obtain values for $M$ which are as tight as possible.

\subsection{Combinatorial Disjunctive Constraints}

When the disjunctive constraint \eqref{eq:general disjunction} is a union of $\mathscr{V}$-polyhedra (e.g., Fig. 1), it suffices to consider only
their extreme points.
Let $J = \bigcup_{i=1}^d \text{ext}(P^i)$ denote the \emph{ground set} and take $\mathscr{S} = \{ \text{ext}(P^i) \}_{i=1}^d \subseteq 2^J$ as the collection of extreme points for each of the polyhedra. 

\begin{figure}[thpb]
  \centering
  \includegraphics[width=0.55\linewidth]{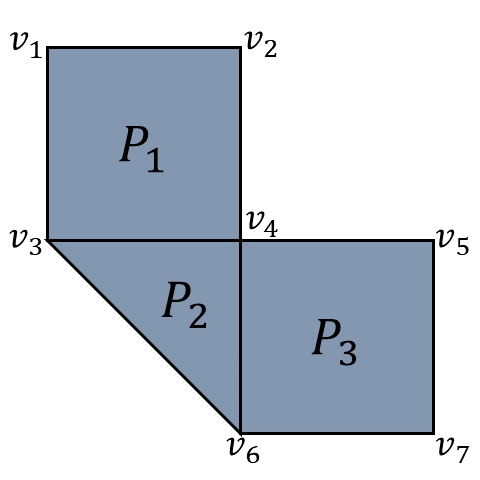}
  \caption{A nonconvex region $\Omega$ defined by a union of polytopes.}
  \label{fig:omega partition with labels}
\end{figure}

\begin{definition}[Definition 1, \cite{CDC}]
\label{def:CDC}
A \emph{combinatorial disjunctive constraint} (CDC) induced by the set $\mathscr{S}$ is
\begin{equation} \label{eq:CDC}
\lambda \in \text{CDC}( \mathscr{S} ) \eqdef \bigcup\limits_{S \in \mathscr{S}} Q(S), 
\end{equation}
where 
$Q(S) \eqdef \{ \lambda \in \Delta^{J} : \lambda_{J \backslash S} = 0 \}$ 
is the face that $S \subseteq J$ induces on the standard simplex.
\end{definition}

In other words, Def. \ref{def:CDC} states $\lambda$ must obtain support over some $S \in \mathscr{S}$. A corresponding formulation for \eqref{eq:general disjunction} is thus  
\begin{equation}
    \label{eq:convex combinations}
    x \in \left\{ \sum_{v \in J} \lambda_{v} v : \lambda \in \text{CDC}(\mathscr{S})  \right\}.
\end{equation}

We highlight that one advantage of this approach is that we can construct a single, strong formulation for a given CDC($\mathscr{S}$) and use it for other instances whose combinatorial structures are sufficiently equivalent, e.g., path planning in similar environments (see \cite{CDC} for a sufficient condition).

As with \cite{CDC}, we assume $\mathscr{S}$ is irredundant: there do not exist distinct $S,T \in \mathscr{S}$ such that $S \subseteq T$. We say that a set $S \subseteq J$ is a \emph{feasible set} with respect to CDC($\mathscr{S}$) if $S \subseteq T$ for some $T \in \mathscr{S}$ and that it is an \emph{infeasible set} otherwise. Furthermore, a \emph{minimal infeasible set} does not contain any infeasible set as a proper subset. 

Standard ideal formulations for \eqref{eq:CDC} are $\mathscr{O} (|\mathscr{S}| + \sum_{S \in \mathscr{S}} |S|)$ in size \cite{CDC, JeroslowLowe}. The \emph{independent branching} (IB) scheme \cite{IBoriginal, CDC}, however, when one exists, is a logically equivalent way of expressing a CDC as a series of choices between two alternatives, whose formulation size on the order of the number of dichotomies.

\begin{definition}[Definition 2, \cite{CDC}]
An \emph{independent branching} scheme for CDC($\mathscr{S}$) is given by a family of sets ($L^j, R^j$) (where $L^j, R^j \subseteq J$) for $j \in [t]$, where
\begin{equation}
    \text{CDC}(\mathscr{S}) = \bigcap_{j=1}^{t} \left( Q(L^j) \cup Q(R^j) \right).
\end{equation}
\end{definition}
Such an IB scheme is said to have \emph{depth} $t$ and that each $j \in [t]$ yields a corresponding \emph{level}.



CDC($\mathscr{S}$) is thus equivalently represented by $t$ constraints, each requiring a selection among two alternatives, and satisfies the condition that 
$$
T \subseteq J \text{ is a feasible set} \Longleftrightarrow \forall j \in [t]: \ T \subseteq L^j \text{ or } T \subseteq R^j.
$$

\begin{figure}[t]
  \centering
  \includegraphics[width=0.55\linewidth]{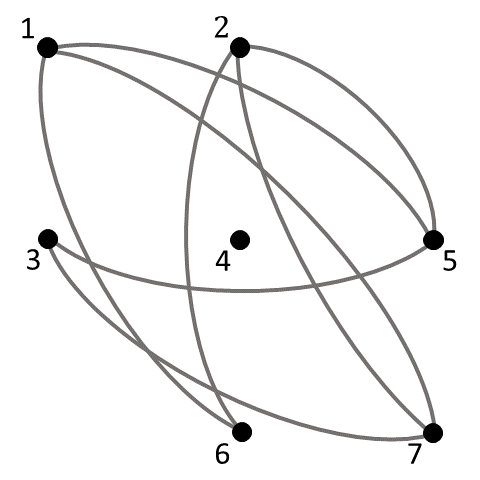}
  \caption{Conflict graph of the polyhedral partition from Fig. \ref{fig:omega partition with labels}. No edge has both endpoints in the same polytope.}
  \label{fig:conflict graph}
\end{figure}

One can obtain IB schemes in an algorithmic manner based on a graphical characterization of CDC($\mathscr{S}$). To illustrate this, we define the \emph{conflict graph} of a CDC($\mathscr{S}$) as 
$G_{\mathscr{S}}^c \eqdef (J, \bar{E}_{\mathscr{S}})$, where $\bar{E}_{\mathscr{S}} \eqdef \{(u,v) \in J^2 : u \neq v, \{u,v\} \text{ is an infeasible set} \}$ \cite{CDC} 
(see Fig. \ref{fig:conflict graph}). Additionally, we require the following notion from graph theory.

\begin{definition}[Definition 4, \cite{CDC}]
    A \emph{biclique cover} of the graph $G = (J,E)$ is a collection of complete bipartite subgraphs $\left\{ G^j \eqdef (A^j \cup B^j, E^j) \right\}_{j=1}^t$ of $G$ that covers all edges of $G$. 
\end{definition}

We will refer to the sets $\{ (A^j, B^j) \}_{j=1}^t$ as a biclique cover.


\begin{proposition}[(Proposition 3, \cite{CDC})]
    If CDC($\mathscr{S}$) is IB-representable and $\{ (A^j, B^j)\}_{j=1}^t$ is a biclique cover for $G_{\mathscr{S}}^c$, then the following is an ideal formulation for CDC($\mathscr{S}$):
    \begin{subequations}
        \label{eq:IBpairwise}
        \begin{align}
            \sum_{v \in A^j} \lambda_v &\leq z_j &\forall j \in [t] \\
            \sum_{v \in B^j} \lambda_v &\leq 1 - z_j &\forall j \in [t] \\
            (\lambda, z) &\in \Delta^J \times \{0,1\}^t.
        \end{align}
    \end{subequations}
\end{proposition}

An intuitive understanding of this formulation is that the construction of infeasible sets is prevented by prohibiting the convex combination of vertices from minimal infeasible sets. 

Consequently, the smallest depth at which an IB scheme can be constructed equates to the minimum biclique cover problem \cite{FishburnBC}. 
While the decision version of this problem is known to be NP-complete \cite{Orlin},
we may instead, however, opt for biclique covers of suboptimal depth, provided they can be obtained efficiently. For example, the trivial biclique cover
$$A^v = \{v\}, \quad B^v = \{u \in J : \{u,v\} \in \bar{E}_{\mathscr{S}} \} \quad \forall v \in J$$ 
allows for a CDC formulation of depth $|J|$ \cite{CDC}. Biclique covers may also be reduced in size by merging bicliques together, provided they satisfy the complete bipartite subgraph definition.

\section{CONTRIBUTIONS}
\label{sec:contributions}

Many obstacle avoidance scenarios naturally give rise to combinatorial disjunctive constraints corresponding to polyhedral partitions in the plane. For this setting, we present: 
\begin{enumerate}
    \item A direct procedure for determining when a CDC admits an IB scheme;
    \item A method for obtaining IB-representable polyhedral partitions;
    \item A polynomial time algorithm for obtaining biclique covers on the class of conflict graphs that arise in these scenarios.
\end{enumerate}



\subsection{IB Scheme Representability}

We begin by defining a polyhedral partition. Consider a (nonconvex) bounded region in the plane $\Omega \subset \mathbb{R}^2$ describing the obstacle-free space an agent may travel in. $\Omega$ can be partitioned into polyhedra $\{P^i \}_{i=1}^d$ such that $\bigcup_{i=1}^d P^i = \Omega$ and whose relative interiors do not overlap (see Fig. \ref{fig:omega partition with labels}). Such a partition will not be unique in general and 
plays a significant role in representability and formulation size.


With a fixed partition, we describe the associated polyhedra in $\mathscr{V}$-form so that the combinatorial disjunctive constraint is characterized by $\mathscr{S} = \{ \text{ext}(P^i) \}_{i=1}^d$ and $J = \bigcup_{S \in \mathscr{S}} S$. As with \cite{CDC}, we also forbid ``internal vertices" by requiring that
\begin{equation}
    \label{eq:internalvertex}
    v \in P^i \Longleftrightarrow v \in \text{ext}(P^i) \quad \forall i \in [d], v \in J .
\end{equation}

\begin{corollary}
    \label{cor:IBscheme plane}
    An IB scheme of a polyhedral partition in the plane exists if and only if each infeasible triplet is not minimal infeasible.
\end{corollary}


One can thus verify IB-representability by enumerating all infeasible triplets and determining if they are minimal infeasible.

Additionally, we present a procedure for partitioning an obstacle-free region $\Omega$ in a manner that is always IB-representable, satisfies the internal vertex condition, and which can be computed quickly without manual user input. A Delaunay triangulation is a triangulation of a given set of vertices (hence, a graph) satisfying certain properties, which we omit for brevity (see \cite{Chew} for a formal definition). A constrained Delaunay triangulation (CDT) requires certain edges be included and, possibly, certain edges within a region be excluded. For the following theorem, we will refer to obstacles with no edges comprising the boundary of $\Omega$ as \emph{internal} obstacles. 

\begin{theorem}
    \label{theorem:CDTpartition}
    Given a set of obstacles with no triangular internal obstacles, the partition of a (nonconvex) obstacle-free region $\Omega$ in the plane yielded from a constrained Delaunay triangulation, which includes the edges corresponding to the boundary of $\Omega$ and prohibits the region enclosed by obstacles from containing diagonal edges, always admits a pairwise IB scheme.
\end{theorem}


CDTs can be computed in $\mathscr{O}(n \log n)$ time (where $n$ is the number of vertices) \cite{Chew}, and software packages such as \emph{Triangle} \cite{shewchuk} have readily available implementations.


\subsection{Obtaining Biclique Covers}

By exploiting the special structure of conflict graphs arising from polyhedral partitions in the plane, we present a polynomial time algorithm for obtaining biclique covers on this class of graphs. We say a graph is \emph{planar} if, intuitively, it can be drawn without any of its edges crossing each other (see \cite{PlanarSep} for a formal definition). The regions enclosed by edges of the graph are referred to as \emph{faces} and can be described in terms of the vertices of the edges that comprise the boundary. A \emph{finite element graph} is any graph formed from a planar embedding of a planar graph by adding all possible diagonal edges to each face, or \emph{element} (i.e., each face becomes a clique). We will relax this definition to allow some faces to remain unaltered.

The following theorem presents a special property of finite element graphs that we leverage in our biclique cover algorithm.

\begin{corollary}[(Corollary 4, \cite{PlanarSep})]
    Let $G$ be any $n$-vertex finite element graph. Suppose no element of $G$ has more than $k$ boundary vertices. The vertices of $G$ can be partitioned into three sets $A,B,C$ such that no edge joins a vertex in $A$ with a vertex in $B$, neither $A$ nor $B$ contains more than $2n/3$ vertices, and $C$ contains no more than $4 \lfloor k/2 \rfloor \sqrt{n}$ vertices.
\end{corollary}

Moreover, Lipton and Tarjan present an algorithm for finding a partition $A,B,C$ in $\mathscr{O}(n)$ time \cite{PlanarSep}.

We denote the conflict graph's complement as $G_{\mathscr{S}}$; it is easy to see that $G_{\mathscr{S}} = \{(u,v) \in J^2 : u \neq v, \{u,v\} \text{ is a feasible set} \}$, and, in particular, $G_{\mathscr{S}}$ will be a finite element graph.
%
%
We can then observe that if we apply the separator algorithm \cite{PlanarSep} on $G_{\mathscr{S}}$ to obtain a partition of its vertices $A,B,C$, the complete bipartite graph induced by $(A,B)$ (assuming neither set is empty) will form a biclique of $G^c_{\mathscr{S}}$.
As there may remain uncovered edges in $G^c_{\mathscr{S}}$, we may likewise apply the separator algorithm to the subgraphs of $G_{\mathscr{S}}$ induced by $A \cup C$ and $B \cup C$, respectively, to obtain further bicliques. 
Applying the separator algorithm to $G_{\mathscr{S}}$ in a divide-and-conquer manner will thus yield a collection of bicliques forming a biclique cover for $G^c_{\mathscr{S}}$.
We provide a high level summary in Algorithm \ref{alg:find BC}.

\begin{algorithm}
    \caption{Biclique cover algorithm}
    \label{alg:find BC}
    \begin{algorithmic}[1]
        \Require $G_{\mathscr{S}}$ a finite element graph 
        
        \State Apply Lipton and Tarjan's separator algorithm for finite element graphs on $G_{\mathscr{S}}$ to obtain a partition $A,B,C$ of the graph's vertices.
        
        \item Apply a postprocessing procedure to ensure $A$ and $B$ are both nonempty (given $G_{\mathscr{S}}$ is not a complete graph). Insert the complete bipartite graph induced by $(A,B)$ into the biclique cover.
        
        \State Recursively apply this algorithm to the subgraphs of $G_{\mathscr{S}}$ induced by the vertices $A \cup C$ and $B \cup C$, respectively, until all subgraphs are complete graphs.
        
    \end{algorithmic}
\end{algorithm}
It is important to note that the planar separator algorithm does not guarantee both $A$ and $B$ be nonempty simultaneously. However, in such a case a postprocessing procedure can be applied to transfer vertices
appropriately.


\begin{theorem}
    \label{theorem:BCconvergence}
    Given a finite element graph $G$, the biclique cover algorithm terminates with a biclique cover for the complement graph $G^c$. Moreover, the complexity of this algorithm is $\mathscr{O}(n^4)$. 
\end{theorem}


\section{COMPUTATIONAL RESULTS}
\label{sec:results}

To demonstrate our framework, we applied our approach to the specific case of footstep planning of a humanoid robot. In these problems, we seek to determine the precise $x,y$ and $\theta$ coordinates of $N$ footsteps en route to a goal pose, such that no step intersects any obstacle and every step is reachable relative to the previous.
We apply much of the formulation techniques proposed by Deits and Tedrake \cite{FootstepDeits}; however, while their implementation leveraged the big-M technique to formulate footstep assignment 
in $xyz\theta$-space, we 
are limited to applying the IB scheme framework in $xy$-space. Nonetheless, we assumed even terrain to forgo $z$ and imposed constraints to limit the change in $\theta$ between steps. 
The final problem is
an MIQCQP 
in which a walking humanoid robot is incentivized to reach a goal destination.

Using our library, \texttt{ClutteredEnvPathOpt}, we conducted experiments comprising of 69 original hypothetical scenarios of 2D obstacles contained within the unit square, which is without loss of generality since our workspace can always be mapped into the unit square. Each scenario contains 3 polyhedral obstacles of various shapes and sizes (which may lie on the boundary), and each obstacle contains at least 4 vertices, with the majority having between 4 and 6 (internal obstacles may not have 3 for IB-representability purposes). We solve each scenario with one, two and three obstacles at a time, as well as with three methods: 1) IB with the original biclique cover computed by our algorithm; 2) IB with the biclique cover obtained after applying a merging procedure to the original; 3) the big-M method (computed as described in Sec. \ref{subsec:bigM} by maximizing over $[0,1]^2$). All scenarios share the same parameters, e.g., 25 footsteps, and begin on the bottom left corner of the unit square, with the goal pose positioned on the top right.

In constructing the optimization problems, we took our obstacle-free space to be the unit square minus the obstacles and use a constrained Delaunay Triangulation to produce an IB-representable polyhedral partition. 
To decrease the number of disjunctions in our footstep assignment constraint, it is possible to merge polytopes by taking their convex hull, provided vertices do not become redundant and the internal vertex condition \eqref{eq:internalvertex} remains satisfied. However, we found that merging the free-space polytopes may result in a partition that is no longer IB-representable. Therefore, we withheld from merging polytopes in our experiments. 

All instances were instantiated with Julia 1.6.1 and solved with Gurobi 9.1.2 on a machine with an Intel Core i7 processor at 1.8 GHz and 8 GB of RAM. For practicality of footstep planning, we set a time limit of 5 minutes. Source code and obstacle files can be found on our library's Github. 

\begin{figure}[t]
    \centering
    \includegraphics[width=0.89\linewidth]{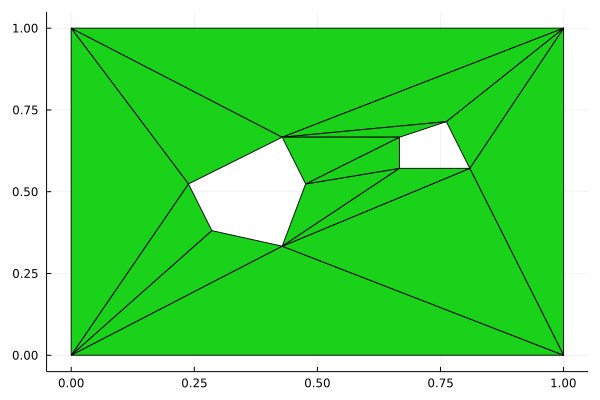}
    \caption{Polyhedral partition from CDT.}
    \label{fig:test free faces}
\end{figure}

\begin{figure}[t]
    \centering
    \includegraphics[width=0.89\linewidth]{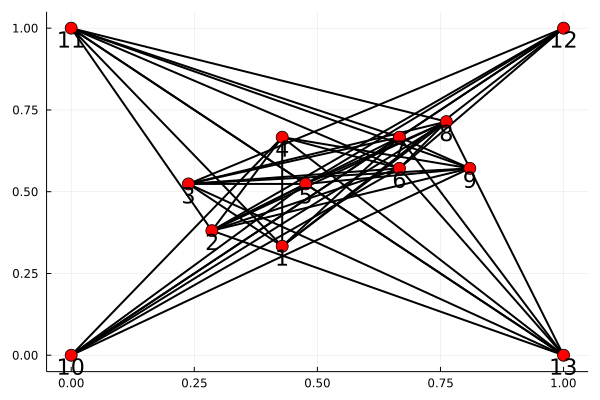}
    \caption{Conflict graph $G^c_{\mathscr{S}}$.}
    \label{fig:test conflict graph}
\end{figure}

\begin{figure}[hbpt]
    \centering
    \includegraphics[width=0.89\linewidth]{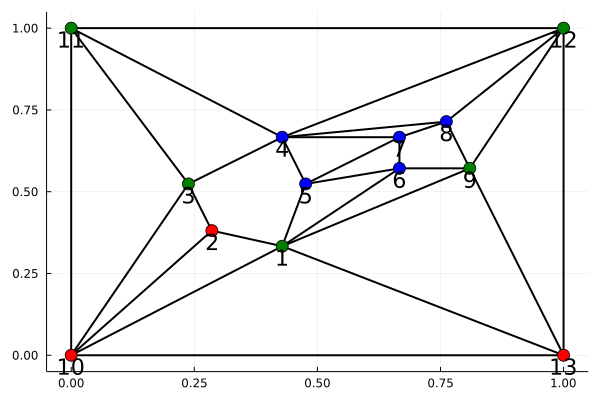}
    \caption{$G_{\mathscr{S}}$. The colored nodes in $G_{\mathscr{S}}$ indicate the vertex partition computed by the planar separator algorithm, with $A$ red, $B$ blue and $C$ green.}
    \label{fig:test FEG}
\end{figure}

\begin{figure}[t]
  \centering
  \includegraphics[width=0.89\linewidth]{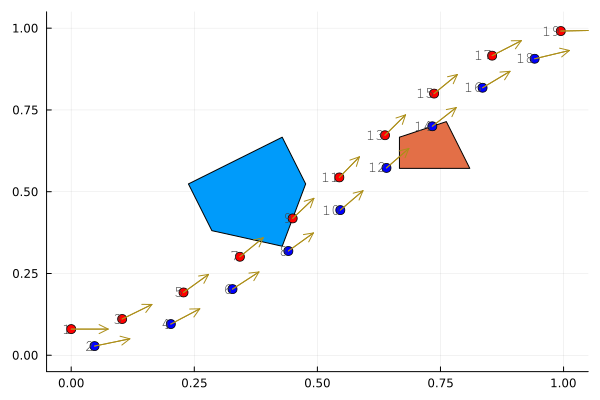}
  \caption{An optimal solution computed via the IB framework.}
  \label{fig:optimal solution}
\end{figure}

To gain a better understanding of our path planning framework, we provide the following example (scenario 47 from our set). The obstacle-free space partition is shown in Fig. \ref{fig:test free faces}, and the associated conflict graph $G^c_{\mathscr{S}}$ and its complement $G_{\mathscr{S}}$ are shown in Figs. \ref{fig:test conflict graph} and \ref{fig:test FEG}, respectively. In this example, there are 13 vertices and 15 ``free'' faces, and the biclique cover obtained by our algorithm is of depth 8 (down from 15 after merging), listed in Table \ref{tab:biclique cover}.
We show an optimal footstep plan for this instance in Fig. \ref{fig:optimal solution}. One can observe the solution opts to trim 6 steps and use 19 from the available 25. One may also notice the path appears to cut the corners of the obstacles, which is consistent with our formulation given we only consider the position of the footsteps when in contact with the ground. This potential shortcoming can be addressed via existing techniques from the literature \cite{CornerCutting}.

\begin{table}[thpb]
    \centering
    \begin{tabular}{ c c c } 
    Level & A & B \\
    \hline
    1 & \{3, 7\} & \{1, 9, 12\} \\
    2 & \{9, 12\} & \{3, 5, 10\} \\
    3 & \{4, 7\} & \{1, 9\} \\
    4 & \{2, 10, 13\} & \{4, 5, 6, 7, 8\} \\
    5 & \{3, 11\} & \{1, 5, 6, 7, 8, 9\} \\
    6 & \{1, 2\} & \{11, 12\} \\
    7 & \{2, 3, 4, 11\} & \{6, 9, 13\} \\
    8 & \{8, 12\} & \{1, 5, 6, 10\} \\
    \hline
    \end{tabular}
    \caption{Biclique cover of conflict graph in Fig. \ref{fig:test conflict graph}.} 
    \label{tab:biclique cover}
\end{table}

In Table \ref{tab:solve results}, we compare IB with the merged cover, IB with the original cover, and the big-M method, across solve times and footstep assignment formulation sizes. We can see that despite having a larger amount of binary variables and general inequalities on average, the big-M method tends to outperform the IB formulations. We highlight that the standard deviations of solve times are larger for the IB approaches as well, reflective of the fact that solve times are competitive in some instances. In fact, our framework was able to solve problems 40-50\% faster in 5 instances and was the fastest in 11 more. 
We can also observe the general trend of increasing solve times with increasing number of obstacles. This is expected: more obstacles result in a larger, more complex conflict graph for the IB case (and likely a larger biclique cover) and more free face inequality systems for the big-M approach. The sizes of these intrinsic properties along with biclique cover reduction statistics are summarized in Table \ref{tab:misc sizes}.

\begin{table*}[t]
    \centering
    \resizebox{\textwidth}{!}{
    \begin{tabular}{ c r r r r r r r r r }
    \hline
    & \multicolumn{3}{c}{1 Obstacle} & \multicolumn{3}{c}{2 Obstacles} & \multicolumn{3}{c}{3 Obstacles} \\
    & IB & IB\_orig & Big-M & IB & IB\_orig & Big-M & IB & IB\_orig & Big-M \\
    \hline
    Fastest & 4 & 4 & 61 & 1 & 3 & 62 & 4 & 0 & 53
    \\
    Timeouts & 14 & 17 & 0 & 49 & 47 & 4 & 55 & 57 & 12
    \\
    Solve Time Avg (s) & 85.06 & 80.71 & 31.18 & 130.36 & 104.12 & 55.74 & 104.18 & 139.26 & 76.02
    \\
    Solve Time Std (s) & 72.11 & 71.41 & 34.54 & 77.18 & 72.22 & 47.22 & 85.12 & 94.27 & 59.1
    \\
    Binary Var. & 115.58 & 179.71 & 208.33 & 192.39 & 351.45 & 359.06 & 241.30 & 492.39 & 496.01
    \\
    Cont. Var. & 214.13 & 214.13 & 0.00 & 317.75 & 317.75 & 0.00 & 415.22 & 415.22 & 0.00
    \\
    Inequalities & 231.16 & 359.42 & 625.00 & 384.78 & 702.90 &  1077.17 & 482.61 & 984.78 & 1488.04
    \\
    \hline
    \end{tabular} }
    \caption{`Fastest' counts and solve time statistics only consider fully solved problems. Formulation sizes are averaged among the 69 test instances.}
    \label{tab:solve results}
\end{table*}

\begin{table}[btph]
    \centering
    \begin{tabular}{ c r r r } 
    \hline
    & 1 Obstacle & 2 Obstacles & 3 Obstacles \\
    \hline
    Vertices & 8.57 & 12.71 & 16.61
    \\
    B.C. Original & 7.19 & 14.06 & 19.70
    \\
    B.C. Merged & 4.62 & 7.69 & 9.65
    \\
    B.C. Reduction (\%) & 35.00 & 44.75 & 50.73
    \\
    Free Faces (F.F.) & 8.33 & 14.36 & 19.84
    \\
    F.F. Halfspaces & 25.00 & 43.09 & 59.52 \\
    \hline
    \end{tabular}
    \caption{Test instance averages of intrinsic properties.} 
    \label{tab:misc sizes}
\end{table}

\section{CONCLUSION}
\label{sec:conclusion}

We have seen that although the IB scheme footstep assignment formulations are (individually) ideal and tend to have fewer binary variables and general inequalites than our respective big-M formulation, it is generally outperformed by the latter. There are a number of possible explanations for this observation.


We hypothesize the strongest factor in the performance discrepancy comes from Gurobi's internal optimization procedures, which include sophisticated branching variable selection techniques, cutting planes, presolve, and numerous heuristics, all of which can reduce the size of the branch-and-cut search tree \cite{GurobiWeb, 50years}. 
A good first direction would be to empirically investigate with other MIQCQP solvers.

Initial feasible solutions can also have a significant affect. The ``nearer'' an initial feasible solution is to the optimal solution, the fewer subproblems we can expect in a branch-and-cut procedure 
\cite{50years}. In a small set of experiments conducted to test this hypothesis, we observed shorter solve times for the IB approach than big-M when warm starting with near-optimal solutions computed from a variety of methods.
Moreover, perusing the obstacle scenarios for patterns of performance, we found the IB approaches to be most often competitive (though not always) when the optimal path involved minimal meandering around obstacles, which we hypothesize is due to the feasibility of initial solutions.
A direction of future work is to understand this phenomena experimentally and develop heuristic procedures for obtaining ``good'' warm starting solutions quickly.


Another factor is that the IB scheme introduces auxiliary continuous variables $\lambda$, in addition to the auxiliary binary variables $z$. In contrast, the big-M method only requires auxiliary binary variables. However, given the moderate size of our final formulations and modern solvers' capabilities for solving large-scale optimization problems \cite{50years}, we do not expect this difference in continuous variables to play a significant role. We also highlight our program contains quadratic constraints, for which strong formulation theory may not apply exactly as for the linear case. However, Gurobi generally solves MIQCQPs by transforming them into linear programs and adding separating hyperplanes ad hoc \cite{GurobiQuadratic}. 


Finally, as our approach remains competitive in some instances, we must be careful not to generalize our results to all scenarios, let alone all path planning MIP approaches, without further treatment. It would be beneficial to expand the obstacle scenario test sets, allowing for larger quantities of obstacles, obstacles with more vertices, more diverse start and goal positions, and more diverse configurations of obstacle placements. This would allow us to determine if there are particular settings for which our framework is favorable; for example, in scenarios with a large quantity of obstacles, the size of the biclique covers could remain significantly smaller than the number of free faces. 

\addtolength{\textheight}{-0.75cm}   



\section*{APPENDIX}

We refer the reader to Garcia \cite{GarciaMasters} for proofs to Corollary \ref{cor:IBscheme plane}, Theorem \ref{theorem:CDTpartition} and Theorem \ref{theorem:BCconvergence}.

\section*{ACKNOWLEDGMENTS}

The authors would like to thank Miles Olson for his software contributions to \texttt{ClutteredEnvPathOpt}. The authors would also like to thank Dr. Lydia E. Kavraki and Carlos Quintero Pena for their robotics related input.



\bibliographystyle{IEEEtran}
\bibliography{IEEEabrv,IEEEexample}

\end{document}